\input amstex
\input prepictex
\input pictex
\input postpictex
\magnification=\magstep1
\documentstyle{amsppt}
\TagsOnRight
\hsize=5.1in                                                  
\vsize=7.8in
\define\R{{\bold R}}

\vskip 2cm
\topmatter

\title Signature via Novikov numbers \endtitle

\rightheadtext{Signature via Novikov numbers}
\leftheadtext{Michael Farber}
\author  Michael Farber\endauthor
\address
School of Mathematical Sciences,
Tel-Aviv University,
Ramat-Aviv 69978, Israel
\endaddress
\email farber\@math.tau.ac.il
\endemail
\thanks{The research was supported by EPSRC Visiting Fellowship.} 
\endthanks
\abstract{It is shown that the signature of a manifold with a symplectic circle action, having 
only isolated fixed points, equals the alternating sum of the Novikov numbers corresponding to
the cohomology class of the generalized moment map. The same is true for more general fixed point
sets}
\endabstract
\endtopmatter
%-------------------------------------------------------------------------
%----------------------------------------------------------------------

\define\C{{\bold C}}

%\define\cl{{\frak Cl}}

%\define\Ker{\operatorname{Ker}}

%\define\ln{\operatorname{ln}}

\def\<{\langle}
\def\>{\rangle}

\define\pd#1#2{\dfrac{\partial#1}{\partial#2}}

\def\part{\partial}

\NoBlackBoxes

\define\ind{\operatorname{ind}}

\documentstyle{amsppt}

\nopagenumbers

\proclaim{1. Theorem} Let $M^{2n}$ be a symplectic manifold with a symplectic circle action
having only isolated fixed points. Then the signature of $M$ is given by
$$\sigma(M) \, =\, b_0(\xi)-b_2(\xi)+b_4(\xi)-b_6(\xi)+\dots+(-1)^n b_{2n}(\xi),\tag1$$
where $\xi\in H^1(M;\R)$ denotes the cohomology class of
the generalized moment map and $b_i(\xi)$ denotes the corresponding Novikov number.
\endproclaim

This generalizes a theorem proven in \cite{JR}, which 
concerns {\it Hamiltonian} circle actions. In the Hamiltonian case
$\xi=0$ and the Novikov numbers
become the usual Betti numbers $b_i(\xi) = b_i(M)$.

\subheading{2} Let us explain the terms used in the statement of Theorem 1. 

First note that any symplectic manifold has a canonical orientation, and the signature $\sigma(M)$
is understood with respect to this orientation. 

Let $\omega$ denote the symplectic form
of $M$. The $S^1$ action is assumed to be {\it symplectic}, which means that for any $g\in S^1$ holds 
$g^\ast\omega =\omega$. Let $X$ denote the vector field
generating the $S^1$-action. Then 
$$\theta = \iota(X)\omega\tag2$$
is a closed 1-form on $M$, which
is called the {\it generalized moment map}. We consider the De Rham cohomology class 
$\xi=[\theta] \in H^1(M;\R)$ of $\theta$. 

For the definition of the Novikov numbers $b_i(\xi)$ we refer to
\cite{BF}, \cite{F}, \cite{N}.

For $n$ odd both sides of formula (1) vanish. 
Indeed, we
follow the convention that the signature of any $4k+2$-dimensional manifold is zero.
The RHS of (1) vanishes because of 
the relation $b_i(\xi) =b_{2n-i}(\xi)$, which
follows directly from definition 1.2 in \cite{BF} of the Novikov numbers 
and the classical Poincar\'e duality. 

\subheading{3. Proof of Theorem 1} The first part of the proof is identical to the arguments used in
\cite{JR}. The second part uses the Novikov - Morse inequalities established in \cite{BF}.

For any fixed point $p\in M$ of the circle action, we obtain a linear representation of the 
circle $S^1$ on the tangent space $T_p(M)$ having no fixed vectors. This defines a complex structure
on $T_p(M)$ and hence a canonical orientation. We will set $\eta(p)= + 1$ (or $\eta(p) = -1$)
depending whether
this orientation coincides (or is opposite)
with the orientation on $T_p(M)$ given by the symplectic form.

{\it Example:} For the standard rotation of the two-sphere $S^2$ around the North 
and the South poles $p_N$ and $p_S$, holds $\eta(p_N)=1$ and $\eta(p_S)=-1$.

In \cite{JR} it is proven, using the Atiyah - Bott fixed point theorem, that
$$\sigma(M) = \sum \eta(p),\tag3$$
where the sum is taken over all the fixed points $p$ of the circle action.

Let us apply theorem 0.3 of \cite{BF} to the generalized moment map $\theta$ given by (2). The critical
points of $\theta$ are precisely the fixed points $p$ of the circle action.  It is well known (cf. \cite{A},
\S\S 2.1, 2.2 and also \cite{JR}) that all these critical points are non-degenerate. 
Also, the indices of the
critical points of $\theta$ are all even and 
$$\eta(p) \, = \, (-1)^{\ind(p)/2},\tag4$$
cf. \cite{A}, \cite{JR}. Applying Theorem 0.3 of \cite{BF}, we obtain that the Morse counting polynomial
$$\Cal M_{\theta}(\lambda)=\sum c_j(\theta) \lambda^j$$
has only terms of even degree 
(where $c_j(\theta)$ denotes the
number of critical points of $\theta$ having index $j$). Hence the polynomial $\Cal Q(\lambda)$ in
(0.7) of \cite{BF} vanishes (by the Morse lacunary principle). 
Therefore we obtain that all odd-dimensional Novikov numbers 
$b_{2j-1}(\xi)$
vanish and for the even-dimensional Novikov numbers 
$$b_{2j}(\xi)=c_{2j}(\theta).\tag5$$
From (3) and (4) we have $\sigma(M) = \sum (-1)^j c_{2j}(\theta)$. Combining this with
formula (5) completes the proof. \qed

We will consider below another generalization of the theorem of \cite{JR}:
we will allow fixed point sets of symplectic circle actions 
of a more general nature.

An oriented manifold $N$ will be said to be {\it an $i$-manifold} (where $i=\sqrt{-1}$)
if all the odd-dimensional Betti numbers $b_{2j-1}(N)$ vanish and the signature $\sigma(N)$ equals
$\sum_{j=0}^{\dim(N)/2} (-1)^j b_{2j}(N)$. In other words, the signature of $N$ equals evaluation
at $\lambda =i$ of the Poincar\'e polynomial of $N$. 

A point (with its canonical orientation) is an $i$-manifold. Any complex projective space $\C P^n$
(with its complex orientation) is an $i$-manifold. If $N_1$ and $N_2$ are $i$-manifolds then so is
their product $N_1\times N_2$.

The result of \cite{JR} can be restated as follows: {\it
any symplectic manifold admitting a Hamiltonian circle action with isolated fixed points is an 
$i$-manifold.} This gives many examples of $i$-manifolds. 

For more examples of $i$-manifolds we refer to D. Metzler \cite{M}.

Theorem 1 remains true if we allow $i$-manifolds 
(instead of isolated points) as connected components of the circle action:

\proclaim{4. Theorem} Let $M^{2n}$ be a symplectic manifold with a symplectic circle action.
Suppose that each component of the fixed point set is an $i$-manifold. 
Then the signature of $M$ can be expressed in terms of the Novikov numbers as follows
$$\sigma(M) = \sum_{j=0}^n (-1)^j b_{2j}(\xi),\tag6$$
where $\xi\in H^1(M;\R)$ is the cohomology class of the generalized moment map.
\endproclaim
Note that for each connected component $Z\subset M$ of the fixed point set of the
circle action, 
$\omega|_Z$ is a symplectic form on $Z$, cf. \cite{A}, where $\omega$ denotes
the symplectic form on
$M$. Hence $Z$ has a canonical orientation and in Theorem 5
we assume that $Z$ is an $i$-manifold with respect to this orientation.

\subheading{5. Proof of Theorem 4} 
For each connected component $Z\subset M$ of the fixed point set we will define
a sign $\eta(Z)=\pm 1$ as follows. The orientation of the tangent bundle to $M$ and the orientation
of the tangent bundle to $Z$ (determined by the symplectic form) define an orientation of the normal
bundle $\nu(Z)$ to $Z$ in $M$. Another orientation of $\nu(Z)$ is determined by the circle action:
for any point $p\in Z$
the circle $S^1$ acts on $\nu_p(Z)$ with no fixed vectors and hence it defines a complex structure and
an orientation of 
$\nu_p(Z)$. We will set $\eta(Z)=1$ if these two orientations of $\nu_p(Z)$ agree; we will set 
$\eta(Z)=-1$ if these two orientations of $\nu_p(Z)$ are opposite.

Using the arguments similar to \cite{JR} one obtains the following formula (generalizing (3))
for the signature
$$\sigma(M) = \sum_{Z} \eta(Z)\sigma(Z),\tag7$$
where the sum is taken over the set of connected components $Z$ 
of the fixed point set. 
Using Theorem 0.3 of \cite{BF} and the lacunary principle as above we find
$$\sum_Z \lambda^{\ind(Z)}\Cal P_Z(\lambda) = \sum_{j\ge 0}\lambda^j b_j(\xi),\tag8$$
where the sum on the left is taken over the connected components $Z$ of the fixed point set and 
$\Cal P_Z(\lambda)$ denotes the Poincar\'e polynomial of $Z$. Hence we obtain (using our assumptions that $\Cal P_Z(\lambda)$ does not involve odd powers) 
that all odd-dimensional Novikov numbers $b_{2j-1}(\xi)$ vanish.

Similarly to (4) we have
$$\eta(Z) = (-1)^{\ind(Z)/2}.\tag9$$
Hence, substituting $\lambda=i$ in (8) and using our assumption 
$\sigma(Z) = \Cal P_Z(i)$
we obtain (6). \qed

\proclaim{6. Corollary} 
Let $M^{2n}$ be a symplectic manifold with a Hamiltonian circle action such 
that each component of the fixed point set is an $i$-manifold. 
Then $M$ is an $i$-manifold. \qed
\endproclaim

\enddocument